\documentclass{amsart}

\usepackage{amscd,amssymb,amsmath}
\usepackage{mathrsfs}
\usepackage{graphics}
\usepackage[all]{xy}

\newtheorem{theorem}{Theorem}[section]
\newtheorem*{maintheorem}{Theorem}
\newtheorem{lemma}[theorem]{Lemma}
\newtheorem{proposition}[theorem]{Proposition}
\newtheorem{corollary}[theorem]{Corollary}

\theoremstyle{definition}

\newtheorem{remark}[theorem]{Remark}
\newtheorem*{acknowledgement}{Acknowledgement}

\theoremstyle{remark}

\renewcommand{\labelenumi}{(\roman{enumi})}

\DeclareFontFamily{U}{wncy}{}
\DeclareFontShape{U}{wncy}{m}{n}{<->wncyr10}{}
\DeclareSymbolFont{mcy}{U}{wncy}{m}{n}
\DeclareMathSymbol{\Sh}{\mathord}{mcy}{"58}

\newcommand\mylabel[1]{\label{#1}}

\newcommand{\NN}{\mathbb{N}}
\newcommand{\ZZ}{\mathbb{Z}}
\newcommand{\QQ}{\mathbb{Q}}

\newcommand{\CC}{\mathbb{C}}

\newcommand{\FF}{\mathbb{F}}

\newcommand{\PP}{\mathbb{P}}

\newcommand  {\shC}     {\mathcal{C}}

\newcommand  {\shT}     {\mathcal{T}}

\def\stX{\mathscr{X}}

\newcommand  {\Br}      {\operatorname{Br}}

\newcommand  {\BGL} 	{\operatorname{BGL}}

\newcommand  {\BPGL} 	{\operatorname{BPGL}}
\newcommand  {\BPU}     {\operatorname{BPU}}
\newcommand  {\BU}      {\operatorname{BU}}

\newcommand  {\Card}    {\operatorname{Card}}

\newcommand  {\Divisible}{{\rm Divisible}}

\newcommand  {\Ext}     {\operatorname{Ext}}

\newcommand  {\Gal}     {\operatorname{Gal}}
\newcommand  {\GL}      {\operatorname{GL}}

\newcommand  {\Hom}     {\operatorname{Hom}}

\newcommand  {\dirlim}  {\varinjlim}
\newcommand  {\invlim}  {\varprojlim}

\newcommand  {\lra}     {\longrightarrow}

\renewcommand{\O}       {\mathcal{O}}

\newcommand  {\ord}     {\operatorname{ord}}

\newcommand  {\Pext}    {\operatorname{Pext}}

\newcommand  {\PGL}     {\operatorname{PGL}}

\newcommand  {\PU}      {\operatorname{PU}}
\newcommand  {\PSL}     {\operatorname{PSL}}

\newcommand  {\quadand} {\quad\text{and}\quad}

\newcommand  {\ra}      {\rightarrow}

\newcommand  {\SL}      {\operatorname{SL}}

\newcommand  {\SU}      {\operatorname{SU}}

\newcommand  {\Torsion} {\operatorname{Torsion}}

\def\mydate{\number\day\space\ifcase\month \or January\or February\or March\or 
April\or May\or June\or July\or
August\or September\or October\or November\or December\fi \space\number\year}

\DeclareFontFamily{U}{wncy}{}
\DeclareFontShape{U}{wncy}{m}{n}{<->wncyr10}{}
\DeclareSymbolFont{mcy}{U}{wncy}{m}{n}
\DeclareMathSymbol{\Sh}{\mathord}{mcy}{"58}

\begin{document}

\title[Brauer groups and phantom cohomology]
      {Infinite CW-complexes, Brauer groups and phantom cohomology}

\author[Jens Hornbostel]{Jens Hornbostel}
\address{Fachbereich C, Bergische Universit\"at Wuppertal,
Gaussstra\ss e 20, 42119 Wuppertal}
\curraddr{}
\email{hornbostel@math.uni-wuppertal.de}

\author[Stefan Schr\"oer]{Stefan Schr\"oer}
\address{Mathematisches Institut, Heinrich-Heine-Universit\"at,
40204 D\"usseldorf, Germany}
\curraddr{}
\email{schroeer@math.uni-duesseldorf.de}


\dedicatory{Revised version, 30 August 2013}

\begin{abstract}
Expanding a result of Serre on finite CW-complexes, we show that the Brauer group coincides
with the cohomological Brauer group for arbitrary compact spaces.
Using results from the  homotopy theory of classifying spaces for Lie groups, 
we give another proof of the result of Antieau and Williams that equality does
not hold for
Eilenberg--MacLane spaces of type $K(\ZZ/n\ZZ,2)$.
Employing a result of Dwyer and Zabrodsky, we show the same for 
the classifying spaces
$BG$ where $G$ is an infinite-dimensional $\FF_p$-vector space.
In this context, we also give a formula expressing
phantom cohomology in terms of homology.
\end{abstract}

\maketitle
\tableofcontents
\renewcommand{\labelenumi}{(\roman{enumi})}

\section*{Introduction}
Generalizing and unifying classical constructions, Grothendieck \cite{GB}
introduced the \emph{Brauer group} $\Br(X)$ and the \emph{cohomological Brauer group} $\Br'(X)$
for arbitrary ringed spaces and ringed topoi. Roughly speaking, the elements in the former
are equivalence classes of geometric objects, which can be regarded, among other things, as 
$\PGL_n$-bundles. In contrast, elements in the latter are cohomology classes
 in degree two of finite order, with coefficients in the multiplicative sheaf $\O_X^\times$ of  units  in the
structure sheaf. The machinery of nonabelian cohomology yields an inclusion  
$$
\Br(X)\subset\Br'(X),
$$
and Grothendieck raised the question under which circumstances
this inclusion is an equality. 
This question is particularly challenging in algebraic geometry, where
one works with  the \'etale site of a scheme $X$.
According to  an unpublished result of Gabber, for which de Jong \cite{de Jong 2006} gave an independent proof,
equality $\Br(X)=\Br'(X)$ holds for quasiprojective schemes.
For many schemes, for example  smooth threefolds without ample sheaves, the question
is regarded as wide open. Note that there are "trivial" counterexamples based on nonseparated
schemes (\cite{Edidin; Hassett; Kresch; Vistoli 1999}, Corollary 3.11)
 and that equality holds for normal algebraic surfaces \cite{Schroeer 2001} and    smooth complex surfaces \cite{Schroeer 2005}.

The goal of this paper is to study Grothendieck's question
in a purely topological settings, where $X$ is a CW-complex, endowed with the sheaf of continuous 
complex-valued functions. In this situation, the cohomological Brauer group $\Br'(X)$ can be identified
with the torsion part of 
$$
\Ext^1(H_2(X)/{\rm Divisible},\ZZ).
$$
According to a result of Serre outlined in \cite{GB}, equality
$\Br(X)=\Br'(X)$ holds for \emph{finite  CW-complexes}.
Slightly expanding Serre's result, we show:

\begin{maintheorem}
For each compact space $X$, we have $\Br(X)=\Br'(X)$.
\end{maintheorem}

A construction of B\"odigheimer \cite{Boedigheimer 1983} 
involving ``long spheres'' then implies that for any torsion group $T$,
there is indeed a compact space---usually not admitting a CW-structure---whose Brauer group is $T$.

The main part of this paper, however, is concerned with 
infinite CW-complexes. Then equality between Brauer group and
cohomological Brauer group does not necessarily hold:

\begin{maintheorem}[Antieau and Williams]
Let $X$ be an Eilenberg--MacLane space of type $K(\ZZ/n\ZZ,2)$.
Then the the Brauer group $\Br(X)$ vanishes, whereas the cohomological Brauer group $\Br'(X)$ is cyclic of order $n$.
\end{maintheorem}

Antieau and Williams  \cite{Antieau; Williams 2011a} used multiplicative properties
of the cohomology ring $H^*(\PU(n),\ZZ)$ with respect to the torsion subgroup.
The case $n=2$ was already considered by Atiyah and Segal \cite{Atiyah; Segal 2004}.
Unaware of these results, we had found another proof;
after putting the first version of this paper onto the  arXiv,   Antieau and Williams 
informed us about \cite{Antieau; Williams 2011a}.

Our approach is based  on a fact from the homotopy theory
of classifying spaces of connected simple Lie groups, stating that
any self map $BG\ra BG$ not homotopic to a constant map
induces bijections on rational homology groups. Nontrivial homotopy classes
of selfmaps indeed exists, namely the \emph{unstable Adams operations}
$\psi^k$ first constructed by Sullivan \cite{Sullivan 1970}, and further studied
by Ishiguro \cite{Ishiguro 1987}, Notbohm \cite{Notbohm 1993}, and 
Jackowski, McClure and Olivier \cite{Jackowski; McClure; Oliver 1992}.

We   also show for arbitrary abelian groups $G$,
any Brauer class on the Eilenberg--MacLane space of type $K(G,2)$ must live in the torsion
part of   $\Ext^1(G/{\rm Torsion},\ZZ)$, and such classes vanish on all compact
subsets (``phantom classes''). 
The latter is a special case of the following purely algebraic description
of such phantoms, a kind of Universal Coefficient Theorem, which   might be useful in other contexts:


\begin{maintheorem}
Let $X$ be a CW-complex and $n\geq 0$. Then the subgroup
$$
\Ext^1(H_{n-1}(X)/{\rm Torsion},\ZZ)\subset H^n(X)
$$
comprises precisely those cohomology classes that vanish on all finite subcomplexes.
\end{maintheorem}


It is easy to see that the cohomological Brauer group for Eilenberg--MacLane spaces of
type $K(G,n)$, $n\geq 3$ vanishes.
We  finally analyze the case $n=1$, that is,  classifying spaces $BG=K(G,1)$ for arbitrary discrete groups $G$.
According to  result of Kan and Thurston \cite{Kan; Thurston 1976}, every CW-complex is homotopy equivalent 
to the \emph{plus construction} $(BG)^+ $ for some group $G$, which is usually uncountable, with respect so some
perfect normal subgroup $N\subset G$. In some sense, this reduces Grothendieck's question
on the equality of $\Br(X)\subset\Br'(X)$   to the case $X=BG$. Our third main result is:

\begin{maintheorem}
Let $X=BG$ be the classifying space of an  infinite-dimensional $\FF_p$-vector space $G$.
Then $\Br(X)\subsetneqq\Br'(X)$.
\end{maintheorem}

The proof relies on 
a result of Dwyer and Zabrodsky \cite{Dwyer; Zabrodsky 1987}, who showed that
every bundle over the classifying space of   finite $p$-groups comes from a representation,
together with some facts on unitary representations due to Backhouse and Bradley \cite{Backhouse; Bradley 1972}.
We actually proof a more general version, where  $G$ can be a $p$-primary torsion group
whose basic subgroups $H\subset G$ are infinite. 

The paper is organized as follows:
In Section \ref{topological brauer groups}, we recall some well-known facts
on Brauer groups in the topological context and show that $\Br(X)=\Br'(X)$ holds for compact spaces.
In Section \ref{low dimensions} we give further preparatory material.
Section \ref{eilenberg-maclane}  contains a new proof that
the Brauer group of an Eilenberg--MacLane space of type $K(\ZZ/n\ZZ,2)$ vanishes,
a discussion of Brauer groups for general $K(G,2)$ and a Universal Coefficient Theorem
for ``phantom'' cohomology classes.
In the final Section \ref{plus construction and classifying spaces}, we turn to classifying spaces $BG=K(G,1)$ for discrete groups.
Here the main result is that for many infinite $p$-primary torsion  groups,
for example infinite-dimensional $\FF_p$-vector spaces, the Brauer group of $BG$
is strictly smaller than the cohomological Brauer group.

\begin{acknowledgement}
We whish to thank Benjamin Antieau and Ben Williams for informing us about  the papers \cite{Antieau; Williams 2011a}, \cite{Antieau; Williams 2012}
and \cite{Woodward 1982}, and the referees for suggestions, which helped to improve the paper. 
The second author wishes to thank Wilhelm Singhof for useful discussions.
\end{acknowledgement}

\section{Recollection: Topological Brauer groups}
\mylabel{topological brauer groups}

In this section we  recall and collect some    useful  facts on Brauer groups of topological spaces,
most of which are well-known.
Throughout, we shall encounter both sheaf and singular cohomology;
if not indicated otherwise, cohomology groups are sheaf cohomology groups. 

First  suppose  $X$ is a general ringed topos. 
We write $\PGL_n(\O_X)$ for the sheaf of groups associated to the presheaf $U\mapsto\PGL_n(\Gamma(U,\O_X))$. 
As explained in \cite{GB},  Section 1, 
the obstruction for extending a $\PGL_n(\O_X)$-torsor $\shT$ to a  $\GL_n(\O_X)$-torsor is a cohomology class
$\alpha\in H^2(X,\O_X^\times)$. More precisely, one has an exact sequence 
\begin{equation}
\label{exact sequence}
H^1(X,\O_X^\times)\lra H^1(X,\GL_n(\O_X))\lra H^1(X,\PGL_n(\O_X))\lra H^2(X,\O_X^\times)
\end{equation}
of pointed sets coming from the theory of nonabelian cohomology 
(compare   \cite{Grothendieck 1955}, Chapter V and \cite{Serre 1972}, Chapter I, \S 5). 
We say that $\alpha$ is the {\it obstruction class}
of the torsor $\shT$; it might also be regarded as a   {\it characteristic class}.
Similarly, the obstruction against extending $\shT$ to an $\SL_n(\O_X)$-torsor
is a class $\tilde{\alpha}\in H^2(X,\mu_n(\O_X))$, which maps to $\alpha$, revealing  that $n\cdot\alpha=0$.
Following Grothendieck, we regard the \emph{cohomological Brauer group} $\Br'(X)$
as the torsion part of $H^2(X,\O_X^\times)$. In contrast, the \emph{Brauer group}
$$
\Br(X)\subset\Br'(X)\subset H^2(X,\O_X^\times)
$$ 
is the subset of elements that are  obstruction classes for some
$\PGL_n(\O_X)$-torsors $\shT$ for certain $n\geq 1$.

Another way to see this goes as follows: The $\alpha\in H^2(X,\O_X^\times)$ correspond
to isomorphism classes of \emph{$\O_X^\times$-gerbes} $\stX$ over $X$,   
and $\alpha$ lies in the Brauer group if one may choose $\stX$ as the gerbe of
extensions to $\GL_n(\O_X)$-torsors for some $\PGL_n(\O_X)$-torsor $\shT$
(see \cite{Giraud 1971}, Chapter V, \S 4).
As explained in \cite{de Jong 2006}, this is equivalent to the existence of certain locally free
\emph{twisted sheaf} of rank $n$, where the twisting is with respect to some cocycle representing $\alpha$.

Now let $X$ be a topological space, and assume that the structure sheaf $\O_X$ is the sheaf
of continuous complex-valued functions $\shC_X$. The exponential sequence
$$
0\lra \ZZ\lra\shC_X\lra\shC_X^\times\lra 1,
$$
where the map on the left is $f\mapsto e^{2\pi i f}$, yields an exact sequence
$$
H^2(X,\shC_X)\lra H^2(X,\shC_X^\times)\lra H^3(X,\ZZ)\lra H^3(X,\shC_X).
$$
Recall that a space $X$ is \emph{paracompact} if it is Hausdorff, and every open covering admits
a refinement that is locally finite. For paracompact spaces, the sheaf $\shC_X$ is  soft by the Uryson Lemma, whence
acyclic (see \cite{Godement 1964}, Chapter II, Theorem 4.4.3), thus one has a canonical
identification of $\Br'(X)$ with the torsion part of $H^3(X,\ZZ)$.


Throughout, we are mainly interested in CW-complexes. Then there is a useful interpretation
in terms of singular homology groups   as well. For 
any abelian group $G$, let us write 
$$
\Torsion G\quadand G/\Divisible 
$$
for the torsion subgroup, and the quotient by the maximal divisible subgroup, respectively.

\begin{proposition}
\mylabel{cw complexes}
If $X$ is a CW-complex, then there is a canonical identification
$$
\Br'(X) = \Torsion\Ext^1(H_2(X),\ZZ),
$$
and this equals the torsion part of $\Ext^1(H_2(X)/\Divisible,\ZZ)$ as well.
\end{proposition}

\proof
The first assertion is well-known, and relies on the fact that CW-complexes
are paracompact (Miyazaki's result, see \cite{Fritsch; Piccini 1990}, Theorem 1.3.5),
that sheaf cohomology for the sheaf of locally constant integer-valued functions
coincides with singular cohomology (compare \cite{Bredon 1967}, Chapter III, Section 1),
and the Universal Coefficient Theorem.


From this the second assertion easily follows: let $D\subset H_2(X)$
be the maximal divisible subgroup. Since divisible groups are injective objects, this
is a direct summand, and it remains to check that $\Ext^1(D,\ZZ)$ is torsion free.
Now any divisible group is a direct sum of groups of the form $\QQ$ and 
$\ZZ[p^{-1}]/\ZZ$. But $\Ext^1(\QQ,\ZZ)$ is a $\QQ$-vector space (in fact, of dimension $2^{\aleph_0}$,
compare \cite{Wiegold 1969}), 
and $\Ext^1(\ZZ[p^{-1}]/\ZZ,\ZZ)$ is isomorphic to the group of $p$-adic integers $\ZZ_p$ (see \cite{Weibel 1994}, p.\ 74).
\qed

\medskip
Let us record the following immediate consequence:

\begin{corollary}
\mylabel{vanishing criterion}
Let $X$ be a CW-complex such that $H_2(X)$ is a direct sum of a divisible group
and a free abelian group. Then $\Br'(X)=0$.
\end{corollary}

For a characterization of abelian groups $G$ with torsion free $\Ext^1(G,\ZZ)$, we refer 
to \cite{Eklof; Mekler 2002}, Theorem 2.13. Under suitable suitable set theoretical assumptions, this are precisely
the sums of divisible and free groups \cite{Mekler; Shelah 1993}.

\medskip
Concerning the Brauer group, one can say the following.
First observe that the sheaf of groups $\PGL_n(\shC_X)$ may be regarded as the sheaf of
continuous $\PGL_n$-valued functions. By abuse of notation, we write $\PGL_n$
for the Lie group $\PGL_n(\CC)$ of $n\times n$ invertible matrices modulo scalar matrices.
There is a well-known equivalence of categories between the category of  
locally trivial principal $\PGL_n$-bundles $P\ra X$ and the category of $\PGL_n(\shC_X)$-torsors $\shT$,
sending $P$ to the sheaf of sections $\shT$ (compare \cite{Grothendieck 1955}, Proposition 5.1.1).

In practice, it is sometimes convenient to work with more geometric objects instead.
Let $f:V\ra X$ be a continuous map whose fibers $f^{-1}(x)$, $x\in X$ are homeomorphic to $\CC\PP^{n-1}$.
A \emph{chart} $(U,\psi)$ consists of an open subset $U\subset X$ and an $U$-homeomorphism $\psi:f^{-1}(U)\ra \CC\PP^{n-1}\times U$.
Two charts $(U,\psi)$ and $(U',\psi')$ are \emph{compatible} if the transition maps $\psi'\circ\psi^{-1}$
take values in $\PGL_n$   fiberwise over $U\cap U'$. An \emph{atlas} is a collection of compatible charts $(U_i,\psi_i)$
with $X=\bigcup_i U_i$. A \emph{$\CC\PP^{n-1}$-bundle} in a continuous map $f:V\ra X$  
as above, endowed with a maximal atlas. One may also view them also as \emph{relative Brauer--Severi-varieties} over $X$  
in the sense of Hakim \cite{Hakim 1971}. 
We have an equivalence of categories between the category of locally principal $\PGL_n$-bundles $P\ra X$
and the category of $\CC\PP^{n-1}$-bundles $V\ra X$, sending the locally trivial principal bundle $P$ to the associated locally trivial fiber bundle
$V$ defined as the quotient of $\CC\PP^{n-1}\times P$ by the diagonal left action of $\PGL_n$ 
(compare \cite{Husemoller 1993}, Chapter 4).

Throughout, we shall write $H^1(X,\PGL_n)$ for the set of isomorphism classes of objects $\shT$, $P$, $V$ 
in their respective categories. We simply use the term \emph{$\PGL_n$-bundle}, or just \emph{projective bundle} when it is clear from the
context to which objects we refer. The operation
$\PGL_m\times\PGL_n\ra\PGL_{mn}$
that comes from   tensor product of matrices induces a pairing 
$$
H^1(X,\PGL_m) \times H^1(X,\PGL_n)\ra H^1(X,\PGL_{mn}),
$$
 which 
in turn gives the addition in the Brauer group. We shall write   $P\otimes P'$ and  $V\otimes V'$ 
for the corresponding operations on  bundles, respectively. 

Stabilizing with respect to tensor products with trivial bundles
leads to the colimit $\PGL=\PGL_{\infty}$. Note that this stabilization
is different from the stabilization for vector resp. $\GL_n$-bundles
denoted by $\GL=\GL_{\infty}$ and given by adding trivial bundles.
These different stabilizations are the reason why $\pi_{2k}(\BPGL) \cong \QQ$
is different from $\pi_{2k}(\BGL) \cong \pi_{2k}(\BU) \cong \ZZ$
for all $k \geq 1$ even if  $\pi_{2k}(\BPGL_n) \cong \pi_{2k}(\BGL_n) \cong 
\pi_{2k}(\BU(n))$. Finally, note that tensoring a bundle with a
trivial bundle does not change its obstruction class. 

Given a topological group $G$ admitting a CW-structure, we denote by $BG$ the classifying space for
numerable principal $G$-bundles. We regard it as a connected CW-complex homotopy equivalent
to the Milnor's join construction $G\star G\star\ldots$, together  with a universal bundle,
and write $[X,BG]$ for the pointed set of homotopy classes of continuous maps.

We are mainly interested in the case $G=\PGL_n$ and $G=\PU(n)$. The latter
is the quotient of the unitary group $U(n)$ by the diagonally embedded circle group $U(1)=S^1$,
which is also the quotient of the special unitary group $\SU(n)$ by the diagonally embedded group of $n$-th roots 
of unity $\mu_n=\mu_n(\CC)$. For paracompact spaces $X$ all locally trivial principal bundles are numerable, 
thus 
\begin{equation}
\label{homotopy equivalence}
H^1(X,\PGL_n)=[X,\BPGL_n] = [X,\BPU(n)],
\end{equation}
the latter because $\BPU(n)\ra\BPGL_n$ is a homotopy equivalence.

%
%

\begin{remark}
Some non-numerable locally trivial principal bundles over non-paracompact
Hausdorff spaces are described in \cite{Bredon 1968} and \cite{Schroeer 2013}.
\end{remark}

Let us record the following fact, which is essentially  Serre's result on the Brauer group of finite CW-complexes:

\begin{theorem}
Let $X$ be a compact space. Then $\Br(X)=\Br'(X)$.
\end{theorem}

\proof
FLet $\alpha\in \Br'(X)\subset H^3(X,\ZZ)$ be a cohomology class, say of order $n\geq 1$.
Chose a lift $\tilde{\alpha}\in H^2(X,\ZZ/n\ZZ)$.
On paracompact spaces, sheaf cohomology coincides with \v{C}ech cohomology, according
to \cite{Godement 1964}, Chapter II, Theorem 5.10.1. Moreover, \v{C}ech cohomology on paracompact spaces with countable coefficients can by described via
homotopy classes of continuous maps into Eilenberg--MacLane spaces, according to P.\ Huber's result
\cite{Huber 1961}. Let $Y$ be an Eilenberg--MacLane space of type $K(\ZZ/n\ZZ,2)$,
and $f:X\ra Y$ a continuous map representing $\tilde{\alpha}$.
The image $f(X)\subset Y$ is compact, whence contained in a finite subcomplex $Y'\subset Y$.
By Serre's result \cite{GB}, Theorem 1.6, the restriction of the universal cohomology class
to $Y'$ is the obstruction of some projective bundle $V'\ra Y'$, and
the pullback $f^*(V')$ shows that $\alpha$ lies in the Brauer group.
\qed

\medskip
From this we get an amusing existence result:

\begin{corollary}
For every abelian torsion group $T$,  there is a compact space $X$
so that $\Br(X)$ is  isomorphic to $T$ and that  the inclusion  $\Br(X)\subset\Br'(X)\subset H^3(X,\ZZ)$ are equalities.
\end{corollary}

\proof
Using ``long spheres'', B\"odigheimer  \cite{Boedigheimer 1983}  constructed  compact spaces
realizing preassigned \v{C}ech cohomology groups. In particular,
there is a compact space $X$ with $\check{H}^3(X,\ZZ)\simeq T$.
But on paracompact spaces, \v{C}ech cohomology coincides with sheaf cohomology,
by \cite{Godement 1964}, Chapter II, Theorem 5.10.1, and the result follows.
\qed

\medskip
Note that for infinite abelian torsion groups $T$, there can be  no CW structure on $X$.
For otherwise there are only finitely many cells by compactness, 
and $H^3(X)$ would be finitely generated by the cellular cochain complex.

\section{Some infinite CW-complexes}
\mylabel{low dimensions}

We start by constructing a concrete projective bundle on a particular CW-complex $Y$ of dimension three.
Fix some $n\geq 1$, and let $Y$ be the space obtained by attaching a 3-cell to the
2-sphere along a continuous map $\varphi:S^2\ra S^2$ of degree $n$. We thus have a cocartesian 
diagram 
$$
\begin{CD}
S^2 @>>> D^3\\
@V\varphi VV     @VV\Phi V\\
S^2 @>>> Y
\end{CD}
$$
and regard $Y=e^0\cup e^2\cup e^3$ as a CW-complex with three cells.
The cellular cochain complex (see for example \cite{Whitehead 1978}, Chapter II, Theorem 2.19),
or the cellular chain complex together with the Universal Coefficient Theorem easily yields
$$
\Br'(Y)=  \Torsion H^3(Y) = H^3(Y) = \ZZ/n\ZZ.
$$
By Serre's result, this  group is generated by some $\PGL_{n'}$-bundle. The following argument
shows that one may actually choose $n'=n$:

Let $L$ be a complex line bundle on the 2-skeleton $Y^2=S^2$ with   Chern class $c_1(L)=1$
(that is, the tautological bundle on $\CC\PP ^1$), where the Chern class is regarded as element of  $H^2(S^2)=\ZZ$, 
and consider the complex vector bundle $E=L\oplus\CC^{n-1}$ of rank $n$. Then
$$
\varphi^*(L\oplus\CC^{n-1})= L^{\otimes n}\oplus\CC^{n-1}.
$$
Although the classifying map $Y^2\ra\BGL_n$ for $E$ does not extend to $Y$,
we have $L^{\otimes n}\oplus\CC^{n-1}\simeq L\oplus\ldots\oplus L$.
This is because vector bundles of rank $n$ on the 2-sphere correspond, up to isomorphism,
to elements in $\pi_1(\GL_n)$, the bundles in question have isomorphic determinant,
and the determinant map $\GL_n\ra\CC^\times$ induces a bijection
on fundamental groups.  Since the projectivization of $L\oplus\ldots\oplus L$
is isomorphic to the product bundle $\CC\PP^{n-1}\times S^2$, the classifying map
$Y^2\ra\BPGL_n$ extends to $Y$, which yields a $\PGL_n$-bundle $V\ra Y$.

\begin{proposition}
\mylabel{one bundle}
The obstruction class of the projective bundle $V\ra Y$ constructed above 
generates the group $\Br'(Y)$. 
\end{proposition}

\proof
Fix $0<i<n$. Since $\Br'(Y)$ is cyclic of order $n$, it suffices to check that 
$V^{\otimes i}$ is not the projectivization of a vector bundle.
Suppose to the contrary that it is the projectivization of a vector bundle $E\ra Y$.
Obviously, it has rank $in$, and $\varphi^*(E)$ is trivial.
In light of the exact sequence of pointed sets  (\ref{exact sequence}) applied to $X=Y^2$, we have
$E|Y^2\simeq (L\oplus\CC^{n-1})^{\otimes i} \otimes N$
for some line bundle $N\ra Y^2$. Consequently
$$
\varphi^*(E) = (L^{\otimes n}\oplus\CC^{n-1})^{\otimes i}\otimes N^{\otimes n}.
$$
Taking first Chern classes, we obtain
$
0=c_1\varphi^*(E) = in^i + in^{i+1} c_1(N)$,
thus $n|i$, contradiction.
\qed

\medskip
The following seems to be well-known, but we could not find a reference for it.
One can also deduce it from the preceding construction:

\begin{proposition}
\mylabel{classifying space}
Let $X=\BPGL_n$. Then the cohomological Brauer group $\Br'(X)$ is cyclic of order $n$, the obstruction class
of the universal bundle $P\ra X$ is a generator, and $\Br(X)=\Br'(X)$.
\end{proposition}

\proof
Clearly, the Lie group $\PGL_n$ is connected. As $\PSL_n \to \PGL_n$
is the universal covering with fiber $\mu_n$ the $n$th roots of unity,
its fundamental group is cyclic of order $n$.
Thus the classifying space $X=\BPGL_n$ is simply connected, and Hurewicz shows that
$H_2(X)=\pi_2(X) =\pi_1(\PGL_n) \simeq \ZZ / n \ZZ$. 
Thus the same holds 
for $\Br'(X)$ by Proposition \ref{cw complexes}.
According to Proposition \ref{one bundle},
there is a numerable $\PGL_n$-bundle  over some space whose obstruction class has order $n$.
Thus the same holds for the universal bundle $P\ra X$. If follows that $\Br(X)\subset\Br'(X)$ is 
an equality.
\qed

\medskip
According to a result of Woodward \cite{Woodward 1982}, see also
\cite{Antieau; Williams 2012}, equality $\Br(X)=\Br'(X)$ holds for all CW-complexes of dimension $\leq 4$;
in fact, each  cohomology class $\alpha\in\Br'(X)$ is the obstruction of some
$\PGL_n$-bundle $V\ra X$ with $n=\ord(\alpha)$.
 
This applies, for example, to differential  manifolds
of   dimension at most four, or the underlying topological space $X=S(\CC)$ of
a algebraic $\CC$-scheme or complex spaces $S$ of complex dimension at most two.
Note that all such spaces   admit a CW-structure.
It also applies to the underlying topological space $X=S(\CC)$ of 
complex Stein spaces $S$ of  complex dimension at most four,
because such spaces have the homotopy type of 
a CW-complex of  real dimension at most four by \cite{Hamm 1983}.

In higher dimensions,  the following criterion reduces the problem
to odd-dimen\-sio\-nal CW-complexes:

\begin{lemma}
\mylabel{extension criterion}
Let $X$ be a CW-complex of even dimension $n$, with $(n-1)$-skeleton
$Y=X^{n-1}$. Let $\alpha\in\Br'(X)$ be a class whose restriction $\alpha|Y$
lies in the Brauer group $\Br(Y)\subset\Br'(Y)$. Then $\alpha\in\Br(X)$.
\end{lemma}

\proof
The cases $n=0,2$ are trivial, for then $H_2(X)$ is free, so $\Br'(X)$ vanishes
by Corollary \ref{vanishing criterion}. Now suppose $n\geq 4$.
Choose some $\PGL_d$-bundle $V\ra Y$ whose obstruction class is $\alpha|Y$.
Let $\varphi:\amalg S^{n-1}_\alpha\ra Y$ be the  attaching maps for 
the $n$-cells. The pullbacks $\varphi^*_\alpha(V)$ correspond to
elements in
$$
\pi_{n-2}(\PGL_d)=\pi_{n-2}(\GL_d)=\pi_{n-2}(U(d)).
$$
Replacing $V$ by the tensor product with a suitable  trivial bundle,
we may assume that $2d+1>n-1$. The fibration
$U(d+1)/U(d)=S^{2d+1}$ shows that we are in the stable range.
We than have $\pi_{n-2}(U(d))=\pi_{n-2}(U)=0$, because   $n$ is even, compare \cite{Mimura 1995}, Section 3.1.
Thus the pullbacks $\varphi^*_\alpha(V)$ are trivial,
so the classifying map $Y\ra\BPGL_d$ for $V$ can be extended to $X$.
The resulting   bundle on $X$ has $\alpha$ as obstruction class,
since the restriction map $\Br'(X)\ra\Br'(Y)$ is bijective.
\qed

\begin{proposition}
\mylabel{even cw complexes}
Let $X$ be a finite-dimensional CW-complex that contains no cells of odd dimension $n\geq 5$.
Then $\Br(X)=\Br'(X)$.
\end{proposition}

\proof
Let $X^n\subset X$ be the $n$-skeleton. We check $\Br(X^n)=\Br'(X^n)$ by  induction on $n\geq 4$.
Equality holds for $n=4$ by \cite{Woodward 1982}.
Now suppose $n\geq 5$, and that equality holds for $n-1$. If $n$ is even, we  may apply   Lemma \ref{extension criterion}.
If $n$ is odd, then $X^{n}=X^{n-1}$, and equality holds as well.
\qed

\section{Eilenberg--MacLane spaces}
\mylabel{eilenberg-maclane}

We now examine the Brauer group of Eilenberg--MacLane spaces $X$ of type $K(G,j)$,
which we always regard as connected CW-complexes endowed with a universal cohomology class.
For $j\geq 3$, such spaces have trivial cohomological Brauer group, according
to Corollary \ref{vanishing criterion}. 
In particular, this holds for $K(\ZZ,3)$, although for paracompact $Y$
every torsion cohomology class  in $H^3(Y,\ZZ)$, which correspond to
an element  of the cohomological Brauer group, arises as a pullback of the universal cohomology 
class on $K(\ZZ,3)$.

The cases of interest are $j=1$ and $j=2$. 
By Proposition \ref{cw complexes}, an Eilenberg--MacLane space $X$ of type $K(\QQ/\ZZ,2)$ has trivial cohomological Brauer group,
although any torsion class in degree three on $Y$ arises as the Bockstein of the pullback of
the universal cohomology class. Similarly for the classifying space
$$
\BPGL_\infty = K(\QQ/\ZZ,2)\times K(\QQ,4)\times K(\QQ,6)\times\ldots.
$$

According to  Proposition \ref{cw complexes}, the cohomological Brauer group of $K(\ZZ/n\ZZ,2)$ is cyclic of order $n$.
Antieau and Williams  proved the following result, using   multiplicative properties of
the cohomology ring $H^*(PU(n),\ZZ)$ with respect to the torsion subgroup (\cite{Antieau; Williams 2011a}, Corollary 5.10):

\begin{theorem}
\mylabel{antieau williams}
The Brauer group of $K(\ZZ/n\ZZ,2)$ is trivial.
\end{theorem}

The special case $n=2$ was already considered by Atiyah and Segal, compare  \cite{Atiyah; Segal 2004}, proof of Proposition 2.1 (v).
Another proof relying on the homotopy theory of classifying spaces that might be of independent interest
appears at the end of this section.

What can be said about  Eilenberg--MacLane spaces $X$ of type
$K(G,2)$, where the abelian group $G$ is arbitrary?   Then $H_2(X)=\pi_2(X)=G$, and Proposition \ref{cw complexes} gives an inclusion 
$\Br'(X)\subset\Ext^1(G,\ZZ)$.

\begin{proposition}
\mylabel{strange brauer}
Assumptions as above. With respect to the inclusion $\Br'(X)\subset\Ext^1(G,\ZZ)$,
the Brauer group   is contained in the subgroup 
$$
\Ext^1(G/{\rm Torsion},\ZZ)\subset\Ext^1(G,\ZZ).
$$
\end{proposition}

\proof
Let $T\subset G$ be the torsion subgroup. 
The short exact sequence of abelian groups $0\ra T\ra G\ra G/T\ra 0$ yields an exact sequence 
$$
\Hom(T,\ZZ)\lra\Ext^1(G/T,\ZZ)\lra\Ext^1(G,\ZZ)\lra\Ext^1(T,\ZZ)\lra 0 
$$
(recall that for abelian groups there are no higher Ext groups). The term on the left vanishes, hence the map on the right indeed is injective. Consider
the collection $G_\alpha\subset G$ of all finite cyclic subgroups. This gives another
short exact sequence $0\ra H\ra\bigoplus_\alpha G_\alpha\ra T\ra 0$  for some torsion group $H$.
In turn, we get a commutative diagram with exact rows and columns
$$
\xymatrix{
         &                        &                             & 0\ar[d]\\
0 \ar[r] & \Ext^1(G/T,\ZZ) \ar[r] & \Ext^1(G,\ZZ) \ar[r]\ar[rd] & \Ext^1(T,\ZZ)\ar[r]\ar[d] & 0\\
         &                        &                             & \prod_\alpha\Ext^1(G_\alpha,\ZZ).
}
$$
Seeking a contradiction, we assume that there is a bundle $V\ra X$ whose obstruction class
viewed as an element $\alpha\in\Ext^1(G,\ZZ)$ does not lie in the subgroup $\Ext^1(G/T,\ZZ)$.
By the preceding diagram, there must be an inclusion $\ZZ/n\ZZ\subset G$ so that the $\alpha$
restricts to  a nonzero element in $\Ext^1(\ZZ/n\ZZ,\ZZ)$. Consider the induced map
$$
f:K(\ZZ/n\ZZ,2)\lra K(G,2) = X
$$
of Eilenberg--MacLane spaces. Then $f^*(V)$ is a bundle with nontrivial obstruction class,
in contradiction to Theorem \ref{antieau williams}.
\qed

\medskip
We conjecture that the Brauer groups for arbitrary Eilenberg--MacLane spaces of type $K(G,2)$ vanish.
In any case, its elements are ``phantoms'':

\begin{theorem}
\mylabel{phantom cohomology}
Let $X$ be a CW-complex and $n\geq 0$. Then the subgroup
$$
\Ext^1(H_{n-1}(X)/{\rm Torsion},\ZZ)\subset H^n(X)
$$
comprises precisely those  classes which vanish on all finite subcomplexes.
\end{theorem}

\proof
According to a generalization  (\cite{Araki; Yosimura 1972}, Corollary 12)  
of Milnor's exact sequence (\cite{Milnor 1962}, Lemma 2), 
we have a natural short exact sequence
$$
0\lra{\invlim}^1 H^{n-1}(X_\alpha)\lra H^n(X)\lra \invlim H^n(X_\alpha)\lra 0,
$$
where $X_\alpha\subset X$ is the ordered set of finite subcomplexes. The Universal Coefficient Theorem
gives   short exact sequences
$$
0\lra \Ext^1(H_{n-2}(X_\alpha),\ZZ)\lra H^{n-1}(X_\alpha) \lra \Hom(H_{n-1}(X_\alpha),\ZZ)\lra 0.
$$
The groups $H_{n-2}(X_\alpha)$ are finitely generated, since $X_\alpha$ are finite CW-complexes.
Hence $\Ext^1(H_{n-2}(X_\alpha),\ZZ)$ are non-canonically isomorphic to the torsion part of $H_{n-2}(X_\alpha)$, whence finite.
According to \cite{Jensen 1972}, Corollary 7.2, the higher derived inverse limits of such groups vanish.
Thus the map 
$$
{\invlim}^1 H^{n-1}(X_\alpha)\lra {\invlim}^1 \Hom(H_{n-1}(X_\alpha),\ZZ)
$$
in the long exact sequence for inverse limits is bijective.
Clearly, we have
$$
\Hom(H_{n-1}(X_\alpha),\ZZ)=\Hom(H_{n-1}(X_\alpha)/{\rm Torsion},\ZZ).
$$
Moreover, the canonical map $\dirlim H_{n-1}(X_\alpha)\ra H_{n-1}(X)$ is bijective.
One easily sees that the induced map
$$
\dirlim  (H_{n-1}(X_\alpha)/{\rm Torsion})\ra H_{n-1}(X)/{\rm Torsion}
$$
is bijective as well. By Jensen's observation (\cite{Jensen 1972}, page 37), we have an identification
$$
{\invlim}^1 \Hom(H_{n-1}(X_\alpha)/{\rm Torsion},\ZZ)=\Pext^1(H_{n-1}(X)/{\rm Torsion},\ZZ),
$$
where  the right hand side is the group of isomorphism classes of pure extensions. Recall that an extension 
in the category of abelian groups $0\ra G'\ra G\ra G''\ra 0$ is called \emph{pure} if it
remains exact after tensoring with arbitrary groups. 
The isomorphism classes of pure extensions form a subgroup
$$
\Pext^1(G'',G')\subset\Ext^1(G'',G').
$$
According to \cite{Fuchs 1970}, \S 53.3, 
this subgroup is the first Ulm subgroup.
Recall that for an abelian group $E$, the \emph{first Ulm subgroup} is $U^1(E)=\bigcap_{m\geq 1}(mE)$.
However, for any torsion free abelian group $F$, the group $\Ext^1(F,\ZZ)$ is divisible 
(compare \cite{Fuchs 1970}, page 223, point (I)),
hence the latter coincides with its own first Ulm subgroup. Combining these observations, one infers
the assertion.
\qed

\begin{remark}
The equality  ${\invlim}^1 H^{n-1}(X_\alpha)=\Pext^1(H_{n-1}(X),\ZZ)$ was already shown by M.\ Huber and Meier \cite{Huber; Meier 1978}.
For a general treatment of ``phantoms'', we refer to McGibbon \cite{McGibbon 1995} and  Rudyak  \cite{Rudyak 1998}, Chapter III.
\end{remark}

\medskip
\emph{Another proof for Theorem \ref{antieau williams}:}
Let $X$ be an Eilenberg--MacLane space of type $K(\ZZ/n\ZZ,2)$.
Consider the Bockstein exact sequence
$$
H^2(X,\ZZ)\lra H^2(X,\ZZ/n\ZZ)\lra H^3(X,\ZZ)\stackrel{n}{\lra} H^3(X,\ZZ).
$$
The term on the left vanishes, by Hurewicz and the Universal Coefficient Theorem.
In turn, the Bockstein map $H^2(X,\ZZ/n\ZZ)\ra H^3(X,\ZZ)$ is injective, and its
image equals $\Br'(X)\subset H^3(X,\ZZ)$, which is cyclic of order $n$.  
Denote by $\tilde{\alpha}\in H^2(X,\ZZ/n\ZZ)$   the universal cohomology class, which generates
this cyclic group, and $\alpha\in\Br'(X)$ the corresponding 
cohomological Brauer class. 
Seeking a contradiction, we suppose that  there is a  non-zero multiple $\beta=r\cdot\alpha$ lying in $\Br(X)$. 
Choose a $\PGL_m$-bundle $V\ra X$  whose obstruction class is $\beta$.  
As $m \cdot \beta = 0$, the number $m$ is a multiple of $\ord(\beta)=n/\gcd(n,r)$.
Tensoring $V$ with the trivial $\PGL_{\gcd(n,r)}$-bundle, we  may assume 
that $n|m$.

The existence of this bundle $V$ has the following consequence:
Let $Y$ be an arbitrary paracompact space  and $\gamma\in \Br'(Y)$  an element of order $n$.
Viewing it as a torsion class in $H^3(Y,\ZZ)$, we may choose some $\tilde{\gamma}\in H^2(Y,\ZZ/n\ZZ)$ mapping to $\gamma$
under the Bockstein. This $\tilde{\gamma}$ corresponds to a homotopy class of continuous maps $h:Y\ra X$,
and the $\PGL_m$-bundle $h^*(V)$ on $Y$ has obstruction $r\cdot\gamma$.  
We will reach a contradiction below by exhibiting  a paracompact space $Y$ and 
a  class $\gamma\in\Br'(Y)$ 
of order $n$ so that $r\cdot\gamma$ does not come from a   
$\PGL_m$-bundle.

We construct such a space $Y$ as a relative CW-complex by attaching a single cell to 
the classifying space $B=\BPGL_m$. 
Using that $\pi_i(B)$ are finite groups for $i\leq 3$, we infer with Serre's refined Hurewicz Theorem \cite{Serre 1953}
that the canonical map $\pi_4(B)\ra H_4(B)$ becomes bijective after tensoring with $\QQ$.
Since $\pi_4(B)\simeq \ZZ$, there is a continuous map  $\varphi:S^4\ra B$ so that the induced map $H_4(S^4,\QQ)\ra H_4(B,\QQ)$
is bijective. The cocartesian diagram
$$
\begin{CD}
S^4 @>>> D^5\\
@V\varphi VV @VV\Phi V\\
B @>>> Y
\end{CD}
$$
defines a relative CW-complex $ B\subset Y$. 
Then we have a long exact sequence
$$
H^i(Y,B)\lra H^i(Y)\lra H^i(B)\lra H^{i+1}(Y,B)
$$
of  relative singular cohomology groups with integral coefficients.
The  relative cohomology groups $H^i(Y,B)$ can be computed using
the cellular cochain complex, thus vanish for $i\neq 5$
(see for example \cite{Whitehead 1978}, Chapter II, Theorem 2.19).
In particular, the restriction map $H^3(Y,\ZZ)\ra H^3(B,\ZZ)$ is bijective.

The obstruction class $\gamma_u\in\Br(B)$ of the universal bundle on $B$  has order $m$, by Proposition \ref{classifying space}.
Let  $\gamma\in\Br'(Y)$ be the 
unique class restricting to $m/n\cdot\gamma_u$, which   has order $n$.
As explained above, this yields a certain
$\PGL_m$-bundle $h^*(V)$ on $Y$ with obstruction class $r\cdot\gamma$, which in turn is classified by a continuous map 
$f:Y\ra B=\BPGL_m$. Consider the commutative diagram
$$
\xymatrix{
S^4 \ar[r]\ar[d]_{\varphi}  &  D^5\ar[d]^{\Phi}\ar[dr]\\
B \ar@/_2em/ [rr]_g   \ar[r]       &  Y \ar[r]_f              & B.\\
}
$$
The map $g\circ\varphi$ factors over $D^5$, whence is homotopic to a constant map.
We infer that the map $g_*:H_4(B,\QQ)\ra H_4(B,\QQ)$ is zero.

On the other hand, $g:B\ra B$ is not homotopic to a constant map,
because $r\cdot\gamma|B=rm/n\cdot\gamma_u\neq 0$.  Theorem \ref{jom} below implies
that $g_*:H_4(B,\QQ)\ra H_4(B,\QQ)$ is bijective, contradiction.
\qed

\medskip
We have just used the following result of Jackowski,  McClure and Oliver  \cite{Jackowski; McClure; Oliver 1992}:

\begin{theorem} 
\mylabel{jom}
Let G be a compact, connected simple Lie group.
Then for every nontrivial element $g \in [BG,BG]$, there is a positive 
integer $k$ and an $\alpha \in \operatorname{Out}(G)$ such
that $g = B\alpha \circ\psi^k$ where $\psi^k$ is a so-called
unstable Adams operation, having the property that
$H^{2i}(\psi^k,\QQ)$ is multiplication by $k^i$.
\end{theorem}

Actually, the units in the monoid of homotopy classes $[BG,BG]$ is the group 
of outer automorphisms $\operatorname{Out}(G)$,
and   the corresponding quotient has an embedding
$$
[BG,BG]/[BG,BG]^\times\subset\NN
$$
into the multiplicative monoid of the integers. 
Its image consists of $k=0$ and
all $k>0$ prime to the order of the Weyl group of $G$.
These numbers $k>0$ correspond to the unstable Adams operations $\psi^k$,
which are defined in terms of the Galois action of $\Gal(\QQ)$ on 
Grassmannians. First examples of unstable Adams operations were constructed by Sullivan \cite{Sullivan 1970}.
This result apply in particular for the compact Lie group $G=\PU(n)$. 
Since $\BPU(n)\ra\BPGL_n$ is a homotopy
equivalence  and $\PU(n)$ is a compact, connected  simple 
Lie group, it applies to our situation.

\section{Plus construction and classifying spaces}
\mylabel{plus construction and classifying spaces}

In this final section we study Eilenberg--MacLane spaces of type $K(G,1)$, that is, classifying spaces
for discrete groups. Somewhat surprisingly, this turns out to be,  in some sense,  the universal situation.

Let $X$ be a CW-complex, and $N\subset\pi_1(X)$ be normal subgroup that is \emph{perfect}, i.e.,
$N$ coincides with its commutator subgroup $[N,N]$. The \emph{plus construction} $X\subset X^+$
is a relative CW-complex obtained by attaching certain cells of dimension two and three, such that the following holds:
\begin{enumerate}
\item $\pi_1(X)\ra\pi_1(X^+)$ is surjective, with kernel $N$.
\item $H_*(X)\ra H_*(X^+)$ is  bijective.
\end{enumerate}
It has the following universal property: For   every CW-complex  $Y$ and every continuous map $f:X\ra Y$
such that $N$ is contained in the kernel of $\pi_1(X)\ra\pi_1(Y)$, there is a unique homotopy class of continuous maps $f^+:X^+\ra Y$
so that $f^+|_X$ is homotopic to $f$. The plus construction was first introduced by Kervaire \cite{Kervaire 1969},
and   used by Quillen to define the higher $K$-groups for rings \cite{Quillen 1971}.

\begin{proposition}
\mylabel{plus construction}
Let $X$ be a CW-complex. Then the restriction maps $\Br'(X^+)\ra\Br'(X)$ and $\Br(X^+)\ra\Br(X)$
are bijective.
\end{proposition}

\proof
Since $H_2(X)\ra H_2(X^+)$ is  bijective, the induced map on cohomological Brauer groups is bijective,
according to Proposition \ref{cw complexes}.
In turn, the map  on Brauer groups is injective, and it remains to check that every
bundle on $X$ extends to $X^+$.
Since  the classifying space $\BPGL_n$ is simply connected, every continuous mapping $f:X\ra \BPGL_n$
extends to $X^+$  by the universal property of the plus construction.
\qed

\medskip
We thus have canonical identifications
$$
\Br'(X^+)=\Br'(X)\quadand \Br(X^+)=\Br(X).
$$
Now let $Y$ be a connected CW-complex. According to a result of Kan and Thurston \cite{Kan; Thurston 1976},
there is a group $G$, a perfect normal subgroup $N\subset G$, and a homotopy equivalence
$$
(BG)^+\lra Y.
$$
Here $BG$ denotes the classifying space of the discrete group $G$, defined in a \emph{functorial way} as
a geometric realization of the nerve of category with a single object and morphism set $G$,
and $(BG)^+$ is the plus construction with respect to $N\subset G=\pi_1(BG)$.
Note that this group $G$ is often uncountable, 
but in any case we have
$$
\Br'(Y)=\Br(BG)\quadand \Br(Y)=\Br(BG)
$$
by Proposition \ref{plus construction}. Concerning the question of equality  of $\Br(X)\subset\Br'(X)$,
we are thus reduced to the case $X=BG$ for discrete groups $G$. This translates part
of the question from algebraic topology
into group theory, as we have
$$
\Br'(BG)=\Torsion \Ext^1(H_2(G)/{\rm Divisible},\ZZ)
$$
by Proposition \ref{cw complexes}. Recall that   $H_2(G)=H_2(G,\ZZ)$ is often referred to as the 
\emph{Schur multiplier}. For $G$ finite,  the Schur multiplier is finite, and thus non-canonically isomorphic to $\Br'(BG)$.

Let us now recall a fundamental notion from abelian group theory: Suppose  $G$ is an abelian torsion group. A subgroup $H\subset G$ is
called \emph{basic} if the following three conditions hold:

\begin{enumerate}
\item The group $H$ is isomorphic to a direct sum of cyclic groups.
\item The residue class group $G/H$ is divisible.
\item For each integer $n$, the inclusion $nH\subset nG\cap H$ is an equality.
\end{enumerate}

Basic subgroups $H\subset G$ always exists, and the isomorphism class of the abstract group $H$ 
is unique. The theory of basic subgroups goes back to Kulikof \cite{Kulikov 1945},
who treated the case that $G$ is \emph{$p$-primary}, that is, the orders of group elements are $p$-power.
For the general situation, compare \cite{Fuchs 1970}, Chapter VI, Section 33.

\begin{theorem}
\mylabel{classifying space discrete}
Let $p>0$ be a prime, and $G$ be a $p$-primary torsion group whose basic subgroups $H\subset G$ are
infinite. Then the classifying space $X=BG$ has the property $\Br(X)\subsetneqq\Br'(X)$.
\end{theorem}

This applies in particular to the case that $G$ is an infinite-dimensional
$\FF_p$-vector space. The proof relies on a  series of facts, which we  have to establish first.
To begin with, recall that for every abelian group, the addition map $G\times G\ra G$ induces
on the graded homology module $H_*(G)$ the structure of an alternating associative ring.
Therefore, the canonical bijection $G\ra H_1(G)$ yields a ring homomorphism $\Lambda^*(G)\ra H_*(G)$
from the exterior algebra to the homology ring.

\begin{proposition}
\mylabel{group homology}
For every abelian group $G$, the canonical map $\Lambda^2(G)\ra H_2(G)$ is bijective.
If $G$ is torsion and $H\subset G$ is a basic subgroup, the induced map $H_2(H)\ra H_2(G)$
is bijective as well.
\end{proposition}

\proof
The first assertion appears in \cite{Brown 1982}, Chapter V, Theorem 6.4.
Now let $G$ be torsion, and $H\subset G$ a basic subgroup.
The task is to check that the canonical map $$\Lambda^2(H)\lra \Lambda^2(G)$$ is bijective.
We start with surjectivity.
Let $g\wedge g'\in \Lambda^2(G)$. Choose $n>0$ with $ng'=0$. Using that $G/H$ is divisible,
we may write $g=na +h$ with $a\in G$ and $h\in H$. Then $g\wedge g'=h\wedge g'=-g'\wedge h$.
Repeating the argument, we see that $g\wedge g'= h\wedge h'$ for some $h,h'\in H$.

It remains to check injectivity. Suppose we have a nonzero element
$x\in\Lambda^2(H)$. Since $H$ is a direct sum of finite cyclic groups,
so is $\Lambda^2(H)$. This implies that there is an integer $n\neq 0$
with $x\not\in n\Lambda^2(H)$. Since $nH=H\cap nG$, the canonical map of $\ZZ/n\ZZ$-modules
$H/nH\ra G/nG$ is injective (in fact bijective). It follows that the image
of $x$ in 
$$
\Lambda^2(H)\otimes\ZZ/n\ZZ=\Lambda^2(H/nH)=\Lambda^2(G/nG)
$$ 
is nonzero, thus the same holds in $\Lambda^2(G)$.
\qed

\begin{corollary}
\mylabel{classifying brauer}
Let $G$ be an abelian torsion group, with basic subgroup $H$ of the form $H=\bigoplus_{i\in I}\ZZ/n_i\ZZ$.
Then the cohomological Brauer group $\Br'(X)$ of the classifying space $X=BG$
is   isomorphic to the torsion part of 
$\prod_{i<j}\ZZ/(n_i,n_j)$.
\end{corollary}

\proof
We have $H_2(G)=\Lambda^2(H)$ by the Proposition. Using $\Lambda^2(\ZZ/n_i\ZZ)=0$, one obtains  a decomposition
$$
\Lambda^2(H)=\bigoplus_{i<j}\Lambda^1(\ZZ/n_i\ZZ)\otimes\Lambda^1(\ZZ/n_j\ZZ)=\bigoplus_{i<j}\ZZ/(n_i,n_j).
$$
Here we have chosen a total order on the index set $I$. Proposition \ref{cw complexes}, together
with the fact that $\Ext^1(\ZZ/d\ZZ,\ZZ)$, $d\neq 0$ is non-canonically isomorphic to $\ZZ/d\ZZ$, now yields the result.
\qed

\medskip
Next, we need a result that allows us to pass back, in rather special circumstances, from
bundles over classifying spaces to representations. Given a representation $\rho:\pi_1(X)\ra L$
in some Lie group $L$, we obtain an associated $L$-bundle $V\ra X$, which is defined as the
quotient of $L\times X$ modulo the diagonal action of the fundamental group.
Such bundles are called \emph{flat}. We are mainly interest in the case $X=BG$, where $\pi_1(X)=G$,
and $L=\PU(n)\subset\PGL_n$.

\begin{proposition}
\mylabel{classifying flat}
Let $p>0$ be a prime, and $G$ be a group that is an  ascending union
of subgroups $G_0\subset G_1\subset\ldots$. Suppose that the $G_j$ are finite $p$-groups
and that the inclusion homomorphisms $G_j\subset G_{j+1}$
admit  left inverses.
Then every $\PGL_n$-bundle  on the classifying space $X=BG$  is associated
to a representation $G\ra\PU(n)$.
\end{proposition}

\proof
We have $BG=\bigcup_j BG_j$.
Let $V\ra BG$ be a  $\PGL_n$-bundle of rank $n\geq 0$. 
and write $V_j=V|BG_j$ for the restrictions. By a result of Dwyer and Zabrodsky \cite{Dwyer; Zabrodsky 1987},
the canonical maps
$$
\Hom(G_i,\PU(n))\lra [BG_i,\BPU(n)]=[BG_i,\BPGL_n]
$$
are surjective. 
Choose a representation $\rho_j:G_j\ra\PU(n)$ so that $V_j$ is isomorphic to the   bundle
associated to $\rho_j$. For each $j$, the representations
$\rho_{j+1}|G_j$ and $\rho_j$ are conjugate.
This follows from \cite{Whitehead 1978}, Chapter V, Corollary 4.4. Replacing inductively  the $\rho_{j+1}$ 
by suitable $a_j^{-1}\rho_{j+1}a_j$,  we may assume that $\rho_{j+1}|G_j=\rho_j$.
In turn, we get a representation $\rho:G\ra\PU(n)$ with $\rho|G_j=\rho_j$.
Consider the associated  $\PGL_n$-bundle $W\ra BG$.
By construction, the bundles $V,W$ restrict to isomorphic bundles on each $BG_j$.

It remains to verify that two $\PGL_n$-bundles $V,W$ of rank $n$ that
are isomorphic on each $BG_j$ are isomorphic.
Let $V_j,W_j$ be the restrictions to $BG_j$. Since $V=\bigcup V_j$ and $W=\bigcup W_j$,
it suffices to construct compatible bundle isomorphism $f_j:V_j\ra W_j$. We do this by
induction on $j\geq 0$. Choose an arbitrary isomorphism $f_0$, and suppose we already
have a compatible family $f_0,\ldots, f_j$. Let $P_{j+1}\ra BG_{j+1}$ be the $\PGL_n$-torsor
of bundle isomorphism $V_{j+1}\ra W_{j+1}$. Then $P_{j+1}|BG_j$ is isomorphic to $P_j\ra BG_j$. 
By assumption, all these torsors admit a section, thus are isomorphic to the trivial torsor.
To conclude the proof, it   therefore suffices to verify that every continuous function
$BG_j\ra\PGL_n$ extends along the inclusion $BG_j\subset BG_{j+1}$.
But this is trivial, because our left inverses $l_{j+1}:G_{j+1}\ra G_j$ for the inclusions $G_j\subset G_{j+1}$
induce by functoriality continuous maps $Bl_{j+1}:BG_{j+1}\ra BG_j$ that are the identity on $BG_j$.
\qed

\begin{proposition}
\mylabel{finite index}
Let $G$ be an abelian  group, and $V\ra BG$ the $\PGL_n$-bundle associated to a representation $\rho:G\ra\PU(n)$.
Then there is a subgroup $M\subset G$ of finite index so that the restriction
$V|BM$ comes from a vector bundle of rank $n$.
\end{proposition}

\proof
This is essentially a result on unitary representations of discrete groups due to  Backhouse and Bradley \cite{Backhouse; Bradley 1972}. 
For the sake of the reader, we briefly recall
parts of their argument  in our setting.  
For each $g\in G$, choose a matrix $A_g\in \SU(n)$ representing $\rho(g)\in\PU(n)$. Make the choice
so the $A_e$ is the unit matrix for the neutral element $e\in G$. Now define
the \emph{multiplier} $\omega_{g,h}\in \mu_n $ by the formula
$$
A_g\cdot A_h = \omega_{g,h}\cdot A_{gh}.
$$
Then the 2-cochain $\omega_{g,h}$ is a cocycle, and its cohomology class in $H^2(G,\CC^\times)$
is precisely the obstruction against lifting the projective representation 
$\rho$ to a linear representation.

On the other hand, one may use this 2-cocycle to endow the set $\tilde{G}=\mu_\infty\times G$
with the structure of a central extension of $G$, by declaring the group law as 
$$
(\zeta,g)\cdot(\xi,h) =  (\zeta\xi\omega_{g,h},gh).
$$
Here $\mu_\infty\subset\CC^\times$ is the group of all complex roots of unity.
The   representation $\rho:G\ra\PU(n)$ yields a unitary representation
$\tilde{\rho}:\tilde{G}\ra U(n)$, defined via $\tilde{\rho}(\zeta,g)=\zeta A_g$.
We may assume that the latter linear representation is irreducible, by passing
to an irreducible subrepresentation. This replaces $n$ by some $n'\leq n$,
but does not affect the multiplier.  

Using Zorn's Lemma, one checks that there is  a maximal subgroup $M\subset G$ on which
the multiplier is symmetric, that is, $\omega_{g,h}=\omega_{h,g}$ for all $g,h\in M$.
Consider the induced central extension $\tilde{M}=\tilde{G}\times_G M$ of $M$.
Obviously, the group $\tilde{M}$ is abelian, and we may view
$0\ra\mu_\infty\ra \tilde{M}\ra M\ra 0$  as an extension in the category of abelian groups.
The latter extension splits,  because $\mu_\infty$ is divisible and hence $\Ext^1(M,\mu_\infty)$ vanishes.
Consequently, the   representation $\rho:M\ra\PU(n)\subset\PGL_n$ lifts to a linear representation
$M\ra\GL_n$. 

It remains to verify that the subgroup $M\subset G$ has finite index.
This indeed holds by \cite{Backhouse; Bradley 1972}, Theorem 3.
Note that in loc.\ cit.\ the irreducibility of the unitary representation of $\tilde{G}$ enters.
There it  is also shown that all irreducible unitary projective  representation with multiplier $\omega $
have the same dimension $d(\omega)$, and that this dimension coincides with the index of $M\subset G$.  
\qed

\medskip
\emph{Proof of Theorem \ref{classifying space discrete}}: 
Let $p>0$ be a prime and $G$ be a $p$-primary torsion group whose basic subgroups $H\subset G$ are
infinite.
Write $H=\bigoplus_{i\in I}\ZZ/p^{\nu_i}\ZZ$. 
By the proof for Corollary \ref{classifying brauer}, the cohomological Brauer group of $BG$ is the
torsion part of
\begin{equation}
\label{product decomposition}
\Ext^1(\Lambda^2(H),\ZZ) = \Hom(\Lambda^2(H),\QQ/\ZZ)=\prod_{i<j}\Hom(\ZZ/(p^{\nu_i},p^{\nu_j}),\QQ/\ZZ).
\end{equation}
In particular, the canonical map $BH\ra BG$ induces a bijection on cohomological Brauer groups.
Replacing $G$ by $H$, we may assume that $G$ itself is a direct sum of cyclic groups.
Restricting to a direct summand and permuting the summands, we may assume
that the index set $I$ is the set of  natural numbers, and that $\nu_i\leq\nu_j$ for $i\leq j$.

Set $X=BG$, and let $\alpha\in\Br'(X)$ be a torsion element. Write it as  a tuple $\alpha=(\alpha_{ij})_{i<j}$
with respect to the decomposition (\ref{product decomposition}).
For each $i\in I$, let $J_i\subset I$ be the set of $j\in I$ with $i<j$ and  $\alpha_{ij}\neq 0$.
Furthermore, let $I_{\rm fin}\subset I$ be the set of all $i\in I$ with  $J_i$ finite.
\emph{Now suppose that $\alpha$ has the property that $I\smallsetminus I_{\rm fin}$ is finite,
and  that the cardinalities of the sets $J_i$, $i\in I_{\rm fin}$ are  unbounded.}
We claim that any such $\alpha$ does not lie in the Brauer group. 

Suppose to the contrary that
there is a $\PGL_n$-bundle $V\ra X$  whose obstruction class is $\alpha$.
Clearly, the assumptions of  Proposition \ref{classifying flat} hold for
the subgroups $G_j=\bigoplus_{i=0}^j\ZZ/p^{\nu_i}\ZZ$ of $G$; hence our projective bundle is associated to a 
projective representation $\rho:G\ra\PGL_n$.
By Proposition \ref{finite index}, there is a subgroup $M\subset G$ of finite index
so that the projective representation lifts to a linear representation
of the subgroup $M$. Thus $\alpha|\Lambda^2(M)=0$.

To reach a contradiction, it remains to check that $\alpha$ remains nonzero on $\Lambda^2(M)$ for
each subgroup $M\subset G$ of finite index. We do this by induction on the index $[G:M]$.
The case $M=G$ is trivial, so let us assume $M\subsetneqq G$.
Choose a short exact sequence
\begin{equation}
\label{hyperplane}
0\lra G'\lra G\lra \ZZ/p\ZZ\lra 0.
\end{equation}
with $M\subset G'$. We shall see that for some other subgroup $U\subset G'$,
we may apply the induction hypothesis to $M\cap U\subset U$. Indeed, we have
$$
U/(M\cap U)= (U+M)/M \subsetneqq G/M,
$$
so $[U:M\cap U]<[G:M]$, and the problem is to choose $U$ admitting a suitable direct sum
decomposition, so that $\alpha|\Lambda^2(U)$ retains its properties with respect to the new decomposition.

With respect to the given direct sum decomposition $G=\bigoplus_{i\geq 0}\ZZ/p^{\nu_i}\ZZ$, the map on the
right in (\ref{hyperplane}) can be viewed as   a matrix $(\lambda_0,\lambda_1,\ldots)$ with $\lambda_k\in\ZZ/p\ZZ$. 
First, consider the case that $\lambda_k=0$ for almost all $k\in I$, say for all $k\geq k_0$.
Then obviously $U=\bigoplus_{i\geq k_0} \ZZ/p^{\nu_i}\ZZ$ does the job.

Suppose now that $\lambda_k\neq 0$ for infinitely many $k\in I$.
Then there must be such an index with $J_k$ finite. Without restriction, we may assume $\lambda_k=1$.
For $i\neq k$, choose lifts $\lambda_i'\in\ZZ$ of $\lambda_i\in\ZZ/p\ZZ$, and consider $e_i'=e_i-\lambda_i'e_k\in G$.
Also set $e_k'=e_k$. Then the elements $e'_i\in G$, $i\geq 0$ form a new ``basis'', and the $e_i'$, $i\neq k$ generate 
a subgroup $U\subset G$  contained in $G'$.
In $\Lambda^2(G)$, we obviously have
$$
e'_i\wedge e'_j = e_i\wedge e_j - \lambda'_i e_k\wedge e_j + \lambda'_j e_k\wedge e_i.
$$
Decomposing $\alpha =(\alpha_{ij}')_{i<j}$ with respect to the new ``basis'' $e_i'\in G$ as in (\ref{product decomposition}), 
the above formula shows that for all $i> k$, $i\not\in J_k$, the condition $\alpha_{ij}\neq0 $ is equivalent to $\alpha'_{ij}\neq 0$, except  for at most
$\Card(J_k)$ indices $j\in J_k$.
It follows that the new tuple $(\alpha'_{ij})_{i<j}$ has the same properties as 
the original tuple $(\alpha_{ij})_{i<j}$, and the same holds if we restrict to indices $i,j\neq k$.
The induction hypothesis applied to $M\cap U\subset U$ tells us that $\alpha$ remains nonzero on $\Lambda^2(M\cap U)$, thus
in particular on $\Lambda^2(M)$.
\qed



\begin{thebibliography}{ccccc}

%

\bibitem{Antieau; Williams 2011a}
B.\ Antieau, B.\ Williams:
The period-index problem for twisted topological K-theory
Preprint, arXiv:1104.4654.

\bibitem{Antieau; Williams 2012}
B.\ Antieau, B.\ Williams:
On the classification of principal $\operatorname{PU}_2$-bundles over a 6-complex.
Preprint, arXiv:1209.2219.

\bibitem{Araki; Yosimura 1972}
S.\ Araki, Z.-I.\ Yosimura:
A spectral sequence associated with a cohomology theory of infinite CW-complexes.
Osaka J.\ Math.\ 9 (1972), 351--365. 

\bibitem{Atiyah; Segal 2004}
M.\ Atiyah, G.\ Segal:
Twisted K-theory. 
Ukr.\ Math.\ Bull.\ 1 (2004),  291--334.

\bibitem{Backhouse; Bradley 1972}
N.\ Backhouse, C.\ Bradley:
Projective representations of abelian groups.
Trans.\ Amer.\ Math.\ Soc.\ 16 (1972), 260--266.

\bibitem{Boedigheimer 1983}
C.-F.\ B\"odigheimer:
Remark on the realization of cohomology groups. 
Quart.\ J.\ Math. Oxford 34 (1983),   1--5.

\bibitem{Bredon 1967}
G.\ Bredon:
Sheaf theory. 
McGraw-Hill Book Co., New York-Toronto-London, 1967.

\bibitem{Bredon 1968}
G.\ Bredon:
A space for which $H^1(X;Z)\not\approx [X,S^1]$. 
Proc.\ Amer.\ Math.\ Soc.\ 19 (1968), 396--398. 

\bibitem{Brown 1982}
K.\ Brown:
Cohomology of groups. 
Springer, Berlin, 1982.

\bibitem{de Jong 2006}
A.\ de Jong:
A result of Gabber.
Preprint, http://www.math.columbia.edu/$\sim$dejong/

\bibitem{Dwyer; Zabrodsky 1987}
W.\ Dwyer, A.\ Zabrodsky:
Maps between classifying spaces. 
In: J.\ Aguad\'{e} and R.\ Kane (eds.), Algebraic topology, pp.\ 106--119.
Lect.\ Notes  Math.\ 1298.
Springer, Berlin, 1987. 

\bibitem{Edidin; Hassett; Kresch; Vistoli 1999}
D.\ Edidin, B.\ Hassett, A.\ Kresch, A.\ Vistoli:
Brauer groups and quotient stacks. 
Amer.\ J.\ Math.\ 123 (2001), 761--777.

\bibitem{Eklof; Mekler 2002}
P.\  Eklof, A.\ Mekler:
Almost free modules. Set-theoretic methods. 
North-Holland, Amsterdam, 2002.

\bibitem{Fritsch; Piccini 1990}
R.\  Fritsch, R.\ Piccinini:
Cellular structures in topology. 
Cambridge University Press, Cambridge, 1990.

\bibitem{Fuchs 1970}
L.\ Fuchs:
Infinite abelian groups. I.
Academic Press, New York-London, 1970.

\bibitem{Giraud 1971}
J.\ Giraud:
Cohomologie non ab\'elienne. 
Springer, Berlin, 1971.

\bibitem{Godement 1964}
R.\ Godement:
Topologie alg\'ebrique et  th\'eorie des faisceaux.
Hermann, Paris, 1964.

\bibitem{Gray 1966}
B. \ Gray: 
Spaces of the same $n$-type, for all $n$.
Topology 5 (1966), 241--243.

\bibitem{Grothendieck 1955}
A.\ Grothendieck:
A general theory of fibre spaces with structure sheaf.
University of Kansas, Department of Mathematics, Report No.\ 4.

\bibitem{Grothendieck 1957}
A.\ Grothendieck:
Sur quelques points d'alg\`ebre homologique. 
Tohoku Math.\ J.\ 9 (1957), 119--221. 

\bibitem{GB}
A.\ Grothendieck:
Le groupe de Brauer. I. Alg\`ebres d'Azumaya et interpr\'etations diverses.
In: J.\ Giraud (ed.) et al.: Dix expos\'es sur la cohomologie des
sch\'emas, pp.\ 46--189.
North-Holland, Amsterdam, 1968. 

\bibitem{Hakim 1971}
M.\ Hakim:
Topos anneles et schemas r\'elatifs.  
Springer, Berlin-Heidelberg-New York, 1972.

\bibitem{Hamm 1983}
H.\ Hamm:
Zum Homotopietyp Steinscher R\"aume.
J.\ Reine Angew.\ Math.\  338  (1983), 121--135.

%

\bibitem{Huber 1961}
P.\ Huber:
Homotopical cohomology and \v{C}ech cohomology. 
Math.\ Ann.\ 144 (1961), 73--76.

\bibitem{Huber; Meier 1978}
M.\ Huber, W.\ Meier: Cohomology and the derived functors of the inverse limit.
In: M.-A.\ Knus, G.\ Mislin, U.\ Stammbach (eds.), Topology and algebra, pp.\ 155-160.
L'Enseignement Math\'ematique, Geneva, 1978.

\bibitem{Husemoller 1993}
D.\ Husemoller:
Fibre bundles.  
Berlin, Springer, 1993.

\bibitem{Ishiguro 1987}
K.\ Ishiguro:
Unstable Adams operations on classifying spaces. 
Math.\ Proc.\ Cambridge Philos.\ Soc.\ 102 (1987), 71--75.

\bibitem{Jackowski; McClure; Oliver 1992}
S.\ Jackowski, J.\ McClure, B.\ Oliver:
Homotopy classification of self-maps of BG via G-actions. I. 
Ann.\ of Math.\   135 (1992),   183--226.

\bibitem{Jensen 1972}
C.\ Jensen:
Les foncteurs d\'eriv\'es de $\invlim$ et leurs applications en th\'eorie des modules. 
Lect.\ Notes  Math.\   254. Springer, Berlin-New York, 1972.

\bibitem{Kan; Thurston 1976}
D.\ Kan, W.\ Thurston:
Every connected space has the homology of a $K(G,1)$.
Topology 15 (1976),   253--258. 

\bibitem{Kervaire 1969}
M.\ Kervaire:
Smooth homology spheres and their fundamental groups.
Trans.\ Amer.\ Math.\ Soc.\ 144 (1969) 67--72.

\bibitem{Kulikov 1945}
L.\ Kulikov:
On the theory of abelian groups of arbitrary cardinality.
Mat.\ Sb. 16 (1945), 129--162.

\bibitem{McGibbon 1995}
C.\ McGibbon:
Phantom maps. 
In: I.\ James (ed), Handbook of algebraic topology, pp.\ 1209--1257.
North-Holland, Amsterdam, 1995. 

\bibitem{Mekler; Shelah 1993}
A.\ Mekler, S.\ Shelah:
Every coseparable group may be free. 
Israel J.\ Math.\ 81 (1993), 161--178. 

\bibitem{Milnor 1962}
J.\ Milnor:
On axiomatic homology theory.
Pacific J.\ Math.\ 12 (1962), 337--341.

\bibitem{Mimura 1995}
M.\ Mimura:
Homotopy theory of Lie groups. 
In: I.\ James (ed.), Handbook of algebraic topology, pp.\ 951--991.
North-Holland, Amsterdam, 1995.

\bibitem{Notbohm 1993}
D.\ Notbohm:
Maps between classifying spaces and applications.  
J.\ Pure Appl.\ Algebra 89 (1993), 273--294. 

\bibitem{Quillen 1971}
D.\ Quillen:
Cohomology of groups. 
Actes du Congr\`es International des Math\'ematiciens, Tome 2, pp.\ 47--51. 
Gauthier-Villars, Paris, 1971. 

\bibitem{Rudyak 1998}
Y.\ Rudyak:
On Thom spectra, orientability, and cobordism. 
Springer, Berlin, 1998.

\bibitem{Schochet 2003}
C.\ Schochet:
A Pext primer: pure extensions and $\lim^1$ for infinite abelian groups. 
NYJM Monographs, Albany, NY, 2003.

\bibitem{Schroeer 2001}
S.\ Schr\"oer:
There are enough Azumaya algebras on surfaces.
Math.\ Ann.\ 321 (2001), 439--454.

\bibitem{Schroeer 2005}
S.\ Schr\"oer:
Topological methods for complex-analytic Brauer groups.
Topology 44 (2005), 875--894.

\bibitem{Schroeer 2013}
S.\ Schr\"oer:
Pathologies in cohomology of non-paracompact Hausdorff spaces.
Topology Appl.\ 160 (2013), 1809--1815.

\bibitem{Serre 1953}
J.-P.\ Serre:
Groupes d'homotopie et classes de groupes ab�liens. 
Ann. \ of \ Math. \ 58 (1953). 258--294.

\bibitem{Serre 1972}
J.-P.\ Serre:
Cohomologie galoisienne.
Fifth edition. Lect.\ Notes  Math.\ 5. 
Springer, Berlin, 1994.

\bibitem{Sullivan 1970}
D.\  Sullivan:
Geometric topology: localization, periodicity and Galois symmetry. 
Springer, Dordrecht, 2005.

\bibitem{Weibel 1994}
C.\ Weibel:
An introduction to homological algebra. 
Cambridge University Press, Cambridge, 1994.

\bibitem{Whitehead 1978}
G.\ Whitehead:
Elements of homotopy theory. 
Springer-Verlag, New York-Berlin, 1978.

\bibitem{Woodward 1982}
L.\ Woodward:
The classification of principal $\operatorname{PU}_n$-bundles over a 4-complex.
J.\ London Math.\ Soc.\  25 (1982),  513--524. 

\bibitem{Wiegold 1969}
J.\ Wiegold:
Ext(Q,Z) is the additive group of real numbers. 
Bull.\ Austral.\ Math.\ Soc.\ 1 (1969), 341--343.

\end{thebibliography}
\end{document}